\documentclass[12pt]{article}
\usepackage[cp1251]{inputenc}
\usepackage{amssymb,latexsym,amsmath}
\usepackage{amsfonts}
\usepackage[english, russian]{babel}
\usepackage{indentfirst}
\usepackage{graphicx,graphics,verbatim}
\usepackage{cmap}
\usepackage{epigraph}

\newtheorem{theorem}{Теорема}[section]

\newtheorem{utv}{Утверждение}
\newtheorem{dfn}{Определение}

\begin{document}

\title{On two classes of generalized fractional operators\\ (with short historical survey of fractional calculus) }

\date{}

\maketitle

\author{
\begin{center}
\textbf{E.~L. Shishkina}\\
Voronezh State University, Voronezh, Russia\\
ilina\_dico@mail.ru\\
\textbf{S.~M. Sitnik}\\
Belgorod State University, Belgorod, Russia\\
sitnik@bsu.edu.ru\\
\end{center}
}

\newpage

\epigraph{The paper is dedicated to Adam Maremovich Nakhushev --- brilliant man and mathematician}

\selectlanguage{english}

\begin{abstract}
This is a survey paper in two parts. In the first part we list main variants of one--dimensional fractional integrodifferential operators. Also some historical and priority remarks are given. As a special question we consider the impact  of Soviet and Russian researches to fractional calculus and its applications to viscoelasticity theory. We also stress an impact of the Voronezh school of mechanics and viscoelasticity theory, including works of Shermergor, Meschkov, Rossikhin, Shitikova and others.

In the second part of the paper we consider two important special classes of generalized fractional operators, they were thoroughly studied by the authors. We consider Buschman--Erdelyi operators,
this is an important class containing as special cases Riemann--Liouville, Erd\'elyi--Kober, Mehler--Fock operators and some others. They also include classical transmutations for the Bessel differential operator,
namely Sonine and Poisson ones. After that we represent results on another important generalization of
Riemann--Liouville operators, namely fractional powers of the Bessel operators.
\end{abstract}

\newpage

\selectlanguage{russian}

\begin{center}
\textbf{О двух классах операторов обобщенного дробного интегродифференцирования
(с коротким историческим обзором)}
\end{center}

\bigskip

\begin{center}
\textbf{С.~М. Ситник}\\
Белгородский государственный национальный исследовательский университет  (НИУ "БелГУ")\\
sitnik@bsu.edu.ru\\
\textbf{Э. Л. Шишкина}\\
Воронежский университет, Воронеж, Россия\\
ilina\_dico@mail.ru
\end{center}

\bigskip

\epigraph{Работа посвящается памяти Адама Маремовича Нахушева --- замечательного человека\\ и математика}

\bigskip

\begin{center}
АННОТАЦИЯ
\end{center}

\bigskip

Эта статья построена в виде небольшого обзора, соответствующие результаты приводятся без доказательств. В первой части
приведен короткий список различных основных вариантов операторов одномерного дробного интегродифференцирования с краткими историческими указаниями. В том числе уточняются некоторые факты об оригинальных приоритетных работах, взаимосвязях различных классов операторов дробного интегродифференцирования, кратко указывается на различные появившиеся в последнее время фальсификации операторов дробного интегрирования. Отдельно рассмотрен вклад советских и российских учёных в приложения дробного исчисления к задачам вязкоупругости в механике, специально отмечен вклад воронежской школы механиков: Шермергора, Мешкова, Россихина, Шитиковой и их учеников.

Во второй части работы рассматриваются два важных специальных класса обобщённых операторов дробного интегродифференцирования, эти классы были ранее подробно изучены авторами. Сначала приведён список основных результатов по теории операторов Бушмана--Эрдейи, этот класс содержит классические операторы Римана--Лиувилля, Эрдейи--Кобера, преобразование Мелера--Фока и ряд других. Этот класс также содержит  операторы преобразования для сингулярного оператора Бесселя, в том числе операторы преобразования Сонина и Пуассона. Далее приведены  основные результаты для другого важного класса обобщений операторов Римана--Лиувилля, а именно, операторов дробного Бесселевого интегродифференцирования.

\bigskip

\textbf{Ключевые слова:} дробное интегродифференцирование, операторы Римана--Лиувилля, операторы Эрдейи--Кобера, операторы преобразования,  операторы Бушмана--Эрдейи, дробные степени оператора Бесселя

\newpage

\tableofcontents

\section{Введение. Основные формы дробного интегродифференцирования}

\textit{Операторы дробного интегродифференцирования} изучены, например,  в \cite{SKM}, \cite{Nah1}--\cite{Nah3}, \cite{KST}, \cite{Psh1}, \cite{Pod}.

Операторы дробного интегродифференцирования играют важную роль во многих современных разделах математики. Для теории специальных функций важность дробного интегродифференцирования отражена в названии известной статьи \cite{Kir4}: "\textbf{\textit{Все специальные функции получаются дробным интегродифференцированием элементарных функций}}"! (Замечание проф. А.А.\,Килбаса: кроме функций Фокса!). В этом пункте мы приведём список основных операторов дробного интегрирования, а также некоторые краткие исторические указания.

\textit{Дробные интегралы и производные Римана-Лиувилля}
\label{DrI}

\begin{dfn}
	Пусть $\varphi (x)\in L_{1} (a,b)$, тогда интегралы
	\begin{gather}
	 \label{RLI1}
	(I_{a+}^{\alpha } \varphi )(x)=\frac{1}{\Gamma (\alpha )} \int\limits_{a}^{x}\frac{\varphi(t)}{(x-t)^{1-\alpha }}dt, \qquad x>a, \\
	 \label{RLI2}
	(I_{b-}^{\alpha } \varphi )(x)=\frac{1}{\Gamma (\alpha )} \int\limits_{x}^{b}\frac{\varphi(t)}{(t-x)^{1-\alpha }}dt, \qquad x<b,
	\end{gather}
	где $\alpha >0$, называются соответственно \textbf{левосторонним \eqref{RLI1}  и правосторонним \eqref{RLI2} дробными интегралами Римана-Лиувилля} порядка $\alpha$ $(0{<}\alpha)$.
	
	Для функции $f(x)$, $x\in[a,b]$ каждое из выражений
	\begin{gather*}
	(D_{a+}^{\alpha }f)(x)=\frac{1}{\Gamma (n-\alpha )} \left(\frac{d}{dx} \right)^{n} \int\limits _{a}^{x}\frac{f(t)dt}{(x-t)^{\alpha -n+1}}, \\
		(D_{b-}^{\alpha }f)(x)=\frac{1}{\Gamma (n-\alpha )} \left(\frac{d}{dx} \right)^{n} \int\limits _{x}^{b}\frac{f(t)dt}{(t-x)^{\alpha -n+1}},
	\end{gather*}
	где $n=[\alpha ]+1,$ ${\alpha >0}$, называется дробной производной
	Римана-Лиувилля порядка $\alpha$, соответственно левосторонней и
	правосторонней.
\end{dfn}

\begin{dfn}
\textit{Операторы Римана--Лиувилля на полуоси} определяются при $\alpha > 0$ по формулам:
\begin{equation}
\label{161}
\begin{gathered}
		(I_{0+}^{\alpha}f)(x)=\frac{1}{\Gamma(\alpha)}\int\limits_0^x \left( x-t\right)^{\alpha-1}f(t)d\,t,\\
		(I_{-}^{\alpha}f)(x)=\frac{1}{\Gamma(\alpha)}\int\limits_x^\infty \left( t-x\right)^{\alpha-1}f(t)d\,t.
 \end{gathered}
 \end{equation}	
Для остальных значений  $\alpha$ они определяются при помощи аналитического продолжения.

Для функции $f(x)$, $x\in[0,\infty)$ каждое из выражений
\begin{gather*}
(D_{0+}^{\alpha }f)(x)=\frac{1}{\Gamma (n-\alpha )} \left(\frac{d}{dx} \right)^{n} \int\limits_{0}^{x}\frac{f(t)dt}{(x-t)^{\alpha -n+1}}, \\
(D_{-}^{\alpha }f)(x)=\frac{1}{\Gamma (n-\alpha )} \left(\frac{d}{dx} \right)^{n} \int\limits_{x}^{\infty}\frac{f(t)dt}{(t-x)^{\alpha -n+1}},
 \end{gather*}
где $n=[\alpha ]+1,$ ${\alpha >0}$, называется дробной производной
Лиувилля порядка $\alpha$, соответственно левосторонней и
правосторонней.
\end{dfn}

Исторически действительно Лиувилль и Риман первыми рассматривали выражения для определённых выше дробных интегралов, Риман на конечном отрезке, Лиувилль на бесконечном. Однако их рассуждения были нестрогими: Лиувилль необоснованно считал, что любая разумная функция раскладывается в ряд по экспонентам (в современной терминологии в ряд Дирихле), Риман сознательно использовал расходящиеся ряды. По-видимому, первое строгое изложение основ теории дробного интегрирования было дано А.В.Летниковым \cite{Letn}, им же операторы дробного интегрирования были впервые применены в качестве операторов преобразования для решения дифференциальных уравнений обобщённого гипергеометрического типа \cite{Letn, Koor}.

Хорошо известно, что обычное дифференцирование $\frac{d}{dx}$ и
интегрирование $\int\limits_a^x...dt$ являются взаимно обратными
операциями, если дифференцирование применяется слева, т.е.
$\frac{d}{dx}\int\limits_a^x\varphi(t)dt=\varphi(x)$. Однако,
вообще говоря, $\int\limits_a^x\varphi'(t)dt\neq\varphi(x)$ (так
как добавляется постоянная $-\varphi(a)$). Точно так же
$\left(\frac{d}{dx}\right)^nI^n_{a+}\varphi=\varphi$, но
$I^n_{a+}\varphi^{(n)}\neq\varphi$, отличаясь от $\varphi$
многочленом порядка $n-1$. Подобным же образом для дробного
дифференцирования всегда будет
$D^\alpha_{a+}I^\alpha_{a+}\varphi=\varphi$, но
$I^\alpha_{a+}D^\alpha_{a+}\varphi$ не всегда совпадает с
$\varphi(x)$ (так как появляются функции $(x-a)^{\alpha-k}$,
$k=1,2,...,[\alpha]-1$, играющие роль многочленов для дробного
дифференцирования).

\textit{Операторы Эрдейи--Кобера} определяются при $\alpha > 0$ по формулам:
\begin{equation}
\label{162}
\begin{gathered}
(I_{0+;\, 2,y}^{\alpha} f)(x) = \frac{2}{\Gamma(\alpha)}x^{-2(\alpha+y)}
\int\limits_0^x (x^2-t^2)^{\alpha-1}t^{2y+1}f(t)\,dt, \\
(I_{-;\, 2,y}^{\alpha} f)(x) = \frac{2}{\Gamma(\alpha)}x^{2y}
\int\limits_x^{\infty} (t^2-x^2)^{\alpha-1}t^{2(1-\alpha-y)-1}f(t)\,dt,
\end{gathered}
\end{equation}
а при значениях $\alpha > -n$, $n \in \mathbb{N}$ по формулам
\begin{gather*}
(I_{0+;\, 2,y}^{\alpha} f)(x) =x^{-2(\alpha+y)} {\left(\frac{d}{d x^2}\right)}^n
x^{2(\alpha	+y+n)}(I^{\alpha+n}_{0+; \, 2,\, y}f)(x)  \\
(I_{-;\, 2,y}^{\alpha} f)(x) = x^{2y} {\left(-\frac{d}{d x^2}\right)}^n
x^{2(\alpha	-y)}(I^{\alpha+n}_{-; \, 2,\, y-n}f)(x).
\end{gather*}
Для остальных значений  $\alpha$ они определяются при помощи аналитического продолжения, аналогично операторам дробного интегродифференцирования Римана--Лиувилля. Операторы Эрдейи--Кобера играют важную роль в теории операторов преобразования, так как такой вид имеют известные операторы преобразования Сонина и Пуассона, или по уточнённой терминологии монографий  \cite{KS, SS} операторы Сонина--Пуассона--Дельсарта (СПД).

Отметим, что в классической монографии \cite{SKM} случаи выбранных нами пределов интегрирования $0,\infty$ не рассматриваются. В последующей английской версии \cite{SaKiMar} эти особые случаи пределов допускаются, но определения содержат неточности, в частности, приводящие к комплексным величинам под знаком интеграла.

\textit{Дробный интеграл по произвольной функции $g(x)$}:
\begin{equation}
\label{163}
\begin{gathered}
(I_{0+,g}^{\alpha}f)(x)=\frac{1}{\Gamma(\alpha)}\int\limits_0^x \left( g(x)-g(t)\right)^{\alpha-1}g'(t)f(t)d\,t,\\
(I_{-,g}^{\alpha}f)(x)=\frac{1}{\Gamma(\alpha)}\int\limits_x^\infty \left( g(t)-g(x)\right)^{\alpha-1}g'(t)f(t)d\,t,
\end{gathered}
\end{equation}
во всех случаях предполагается, что ${\rm Re}\, \alpha>0$, на оставшиеся значения
$\alpha$ формулы также без труда продолжаются \cite{SKM}.
При этом обычные дробные интегралы (\ref{161}) получаются при выборе в определениях
\eqref{163}  $g(x)=x$, Эрдейи--Кобера (\ref{162}) при наиболее часто рассматриваемом  выборе $g(x)=x^2$  или в общем случае $g(x)=x^p$,
Адамара при $g(x)=\ln x$,  в работе М.М.~Джрбашяна и А.Н.~Нерсесяна \cite{DN58} рассматривается вариант  $g(x)=\exp(-x)$.

А.М.\,Джрбашян обратил наше внимание на тот факт, что операторы дробного интегрирования по функции (\ref{163}) являются частными случаями несколько более общих операторов, которые были введены и изучались его отцом М.М.\,Джрбашяном, см. \cite{DZ1}--\cite{DZ5}, а также монографию \cite{SKM}. В этих работах изучены интегральные представления указанных интегродифференциальных операторов, их обращение и соответствующие интегродифференциальные уравнения с операторами дробного порядка.

Дальнейшие обобщения операторов дробного интегродифференцирования связаны с использованием композиций уже введённых более простых операторов.

\textit{Осреднённый оператор  дробного интегродифференцирования, ассоциированный c оператором $R^t$},  вводится по формуле
\begin{equation}
\label{M1}
I_{MR}^{(a,b)}f=\int\limits_a^b R^t f(t)d\,t,
\end{equation}
где $R^t$ --- произвольный оператор дробного интегродифференцирования порядка $t$.
В случае,  когда $R^t$ является дробным интегралом Римана--Лиувилля, общепринятыми являются названия \textit{континуальный} или \textit{распределённый} дробный интеграл (или производная).
Свойства таких операторов и уравнений с ними подробно изучены в работах А.В.Псху и его учеников, см. \cite{Psh1, Psh2}.

Отметим, что одним из авторов была ранее предложена модификация осреднённого выражения в \eqref{M1}, чтобы она стала совпадать с классическим интегральным средним
\begin{equation*}
\overline{I}_{MR}^{(a,b)}f=\frac{1}{b-a}\int\limits_a^b R^t f(t)d\,t.
\end{equation*}
Этот вариант подправленного определения заинтересовал и был одобрен А.М.Нахушевым, но он не используется.

\textit{Дробная производная А.Н. Герасимова } была введена им \cite{Ger, Nov} по формуле
\begin{equation*}
(D_{G\_}^{\alpha}f)(x)=\frac{1}{\Gamma(\alpha)}\int\limits_x^\infty \left( t-x\right)^{\alpha-1}f'(t)d\,t.
\end{equation*}
Аналогично можно ввести дробные производные Герасимова и для других аналогов дробных операторов Римана--Лиувилля, например, для дробных производных Бесселя, см. недавнюю работу авторов \cite{SS3}.

Очевидная модификация дробной производной Герасимова на случай более высокого порядка приводит к дробным производным и интегралам Герасимова--Капуто

\begin{equation*}
(D_{GC\_}^{\alpha}f)(x)=\frac{1}{\Gamma(\alpha)}\int\limits_x^\infty \left( t-x\right)^{\alpha-1}f^{(n)}(t)d\,t
\end{equation*}
и их вариантам.

Одной из областей, в которой производные Герасимова и Герасимова--Капуто нашли важные применения, является теория вязкоупругости в механике.
Статья А.Н. Герасимова оказалась первой публикацией в  литературе на русском языке, в которой дробная производная использовалась для описания вязкоупругих материалов, кроме того в ней впервые была введена важная дробная производная Герасимова, в этом её важность для теории вязкоупругости.
Отметим книгу О.Г.~Новоженовой \cite{Nov}, где собраны по крупицам биографические сведения и публикации А.Н.~Герасимова, включая пионерскую работу \cite{Ger}.
На приоритет А.Н.~Герасимова впервые внимание математиков обратил А.А.~Килбас в своих лекциях \cite{Kil}. Но для механиков это сделал Ю.Н. Работнов в своей книге "Ползучесть элементов конструкций" (\cite{Rab2}, 1966), где он ссылается на статью Герасимова (\cite{Ger}, 1948).
В более общей международной перспективе в ретроспективной обзорной статье \cite{RS3} приведена подробная таблица, в которую сведены практически все первые работы по использованию дробного исчисления в вязкоупругости, изданные в 1940--70е годы. Из таблицы в этой статье  видно, что последовательность публикаций следующая:  Gemant (1936), Scott Blair (1944), Gross (1947), потом  в 1948 году Герасимов и Работнов (они использовали разные подходы), далее Шермергор (1966, \cite{Sher1}), Мешков (1967, {Mesh2}), Мешков и Россихин (1968--1971). Капуто и Майнарди эти же модели записали позже в 1971 году.
Отметим, что сам М.~Капуто никогда не претендовал на открытие указанной производной, его работы развивали в том числе приложения дробного исчисления в механике.
Называть вопреки исторической справедливости производную его именем предложил его ученик Ф.~Майнарди, и с успехом это название пропагандировал и пропагандирует,
хотя Майнарди знал про приоритет Герасимова с самого начала.

\textit{Дробное дифференцирование Джрбашяна--Нерсесяна}, ассоциированное с последовательностью
\begin{equation*}
\{\gamma_0,\,\gamma_1,\,\dots,\,\gamma_m\}
\end{equation*}
порядка $\sigma$, где
$\sigma=\gamma_0+\gamma_1+...+\gamma_m$,
определяется соотношением
\begin{equation}
\label{comp}
D_{DN}^\sigma=D^{\gamma_0} D^{\gamma_1} \cdots D^{\gamma_m},
\end{equation}
где $D^{\gamma_k}$ --- дробные интеграл и производные Римана--Лиувилля c некоторым началом.
Эти операторы вводились в \cite{DN1}--\cite{DN3}  и затем использовались в \cite{DZ6}--\cite{DZ9}. В первоначальных определениях содержались условия $-1\leq\gamma_0\leq 0,\ \  0\leq\gamma_k\leq 1,\ \  1\leq k \leq m$, при которых в указанных работах исследовались уравнения с такими операторами, но операторы Джрбашяна-- Нерсесяна можно рассматривать при любых параметрах $\gamma_k$ при соответствующем понимании операторов дробного интегродифференцирования Римана--Лиувилля.

Операторы Римана--Лиувилля, Герасимова и Герасимова --Капуто являются специальными случаями операторов Джрбашяна--Нерсесяна.

Операторы Герасимова и Джрбашяна--Нерсесяна послужили образцами для введения в книге Миллера и Росса \cite{MR} более общих \textit{последовательных операторов дробного интегродифференцирования}, в которых композиции в формулах типа \eqref{comp} составлены уже из произвольных операторов дробного интегродифференцирования, см.  обсуждение в книге \cite{Pod}.

Другим путём для обобщений операторов дробного интегродифференцирования является рассмотрение конструкций со всё более сложными специальными функциями в определениях, на этом пути были изучены операторы Сайго \cite{Repin}, Лава и многие другие  \cite{SKM}. По мнению авторов также важными обобщениями дробного интегродифференцирования Римана--Лиувилля являются операторы Бушмана--Эрдейи и дробные операторы Бесселя, которые рассмотрены в данной работе.

Необходимо отметить, что огромный вклад в приложения дробного исчисления внесли советские и российские механики. Здесь следует прежде всего выделить
 значительный вклад академика Ю.Н.~Работнова и его школы, книги и работы которого заложили основы нескольких разделов строительной механики \cite{Rab1,Rab2,RS1}.
На основе его работ, в которых использовались методы дробного исчисления, разрабатывались и были включены в ГОСТы нормативы для строительства фундаментов знаний.
Поэтому без большого преувеличения можно сказать, что все дома и другие здания  советской эпохи 1950--80 годов были построены на фундаменте дробного исчисления.
Существенный вклад в приложения дробного исчисления внесли советские механики Д.~Шермергор, М.~Розовский, А.~Ржаницын, С.~Мешков, Ю.~Россихин, М.~Шитикова  и другие.
О реальной истории приложений дробного исчисления в механике см. обзорные статьи Ю.А.~Россихина и М.В.~Шитиковой \cite{RS2, RS3}.

"Меккой"\, дробного исчисления в России иногда называют город Нальчик, благодаря замечательной школе,  созданной здесь Адамом Маремовичем Нахушевым.
Из его книг многие узнали о существовании дробного исчисления и на этих книгах учились.
А.М.~Нахушевым подготовлено много талантливых учеников, которые продолжают исследования в области дробного исчисления и его приложений.
В настоящее время эту школу возглавляет А.В.~Псху --- также ученик А.М.~Нахушева.

Отдельно отметим вклад механиков Воронежской школы в развитие и приложения дробного исчисления. Наряду с работами С.И.~Мешкова \cite{Mesh1,Mesh2}
существенный вклад внесли Воронежские учёные Ю.А.~Россихин и М.В.~Шитикова, развивая различные разделы механики с использованием идей дробного исчисления. Особенно следует отметить, что известные и широко используемые модели Кельвина--Фойгта и Максвелла с дробными производными в вязкоупругости имеют "Воронежское происхождение"{}, что обычно не отмечается и заменяется неверными приоритетными ссылками. Модель Кельвина--Фойгта с дробными производными в вязкоупругости  была впервые написана  в 1966 г. в статье \cite{Sher1} Т.Д. Шермегора, а модель Максвелла с дробными производными в вязкоупругости была впервые написана  в 1967 г. в статье \cite{Mesh2} С.И. Мешкова, оба они  в это время работали в Воронеже. Далее эти модели изучались в 1968--1971 гг. в работах С.И. Мешкова и Ю.А. Россихина, см. подробнее \cite{RS3}. Только потом Капуто и Майнарди в точности эти же модели записали в 1971 году.

В настоящий момент эта работа продолжается в рамках созданного в Воронежском строительном (ныне политехническом) университете Международного научного центра по фундаментальным
исследованиям в области естественных и строительных наук имени Заслуженного деятеля науки РФ
профессора Россихина Ю.А. Возглавляет центр его директор --- М.В.~Шитикова, идеи и приложения дробного исчисления по-прежнему остаются основными в тематике работ этого центра,
в котором обучаются иностранные студенты и аспиранты из многих стран мира.

Подводя определённые итоги развития дробного исчисления и его многочисленных приложений, в 2019 году была издана в определённом смысле энциклопедия в восьми томах "Handbook of
Fractional Calculus with Applications" \cite{Frac}\!, редакторами томов и авторами основных статей которой выступили известные специалисты в области дробного исчисления и
его приложений: Ю.~Лучко, А.~Кочубей,  В.~Кирякова, Р.~Хилфер, Р.~Горенфло, Ф.~Майнарди, С.~Рогозин, М.~Маламуд, К.~Дитхельм, М.~Шитикова и другие.

Тематика дробного исчисления и его приложения тесно связаны с другим важным разделом --- теорией операторов преобразования. Многочисленные примеры использования дробного исчисления в теории операторов преобразования приведены в монографиях авторов \cite{KS,SS,ShSi,Shi}, сборнике статей под редакцией В.В. Кравченко и С.М.  Ситника \cite{KrSi1}, а также в  \cite{Sita4,KrSi2}.

Вместе с тем стоит отметить, что в последнее время недобросовестными и непрофессиональными авторами вводятся многочисленные обобщения и аналоги дробных производных,
которые не могут считаться полноценными объектами дробного исчисления. Некоторые из них не удовлетворяют полугрупповому свойству  или сводятся к простому умножению
на функцию, другие вводятся по аналогии с  формулой для резольвенты интеграла Римана-Лиувилля и также не имеют основных свойств дробного интеграла или производных.
Выяснению важного вопроса, какие свойства необходимы для вводимых обобщений дробных интегралов и производных, а также критике неудачных вариантов посвящён недавний проект Юрия Лучко \cite{Luch}.

\section{Операторы Бушмана--Эрдейи}

В этом пункте мы изложим основные определения и результаты для операторов Бушмана--Эрдейи различных типов.
Это класс операторов, который при определённом выборе параметров является одновременным обобщением операторов преобразования Сонина и Пуассона, операторов дробного интегродифференцирования Римана--Лиувилля
и Эрдейи--Кобера, а также интегральных преобразований Мелера--Фока. Подробное изложение можно найти в монографиях авторов \cite{KS,SS,ShSi}, обзорах \cite{Sita4,Sita6}, сборнике статей \cite{KrSi1} и статьях \cite{Sita3,Sita5,Sita6,Sita8,Sita9,KrSi2}.

\begin{dfn}
	Операторами Бушмана--Эрдейи первого рода называются интегральные операторы
	\begin{gather}
	\label{2BE1}
	(B_{0+}^{\nu,\mu}f)(x)=\int\limits_0^x \left( x^2-t^2\right)^{-\frac{\mu}{2}}P_\nu^\mu \left(\frac{x}{t}\right)f(t)d\,t,\\
	(E_{0+}^{\nu,\mu}f)(x)=\int\limits_0^x \left( x^2-t^2\right)^{-\frac{\mu}{2}}\mathbb{P}_\nu^\mu \left(\frac{t}{x}\right)f(t)d\,t,\\
	(B_{-}^{\nu,\mu}f)(x)=\int\limits_x^\infty \left( t^2-x^2\right)^{-\frac{\mu}{2}}P_\nu^\mu \left(\frac{t}{x}\right)f(t)d\,t,\\
	\label{2BE2}
	(E_{-}^{\nu,\mu}f)(x)=\int\limits_x^\infty \left( t^2-x^2\right)^{-\frac{\mu}{2}}\mathbb{P}_\nu^\mu \left(\frac{x}{t}\right)f(t)d\,t.\\ \nonumber
\end{gather}
\end{dfn}
Здесь $P_\nu^\mu(z)$---функция Лежандра первого рода \cite{BE1}, $\mathbb{P}_\nu^\mu(z)$---та же функция на разрезе $-1\leq  t \leq 1$, $f(x)$---локально суммируемая функция,
удовлетворяющая некоторым ограничениям на рост при $x\to 0,x\to\infty$. Параметры $\mu,\nu$---комплексные числа, ${\rm Re}\,  \mu <1$, можно ограничиться значениями ${\rm Re}\,  \nu \geq -1/2$.

Интегральные операторы указанного вида с функциями Лежандра в ядрах впервые встретились в работах E.T.~Copson по уравнению Эйлера-Пуассона-Дарбу в конце 1950-х годов.

\textbf{Лемма Копсона.}
Рассмотрим дифференциальное уравнение в частных производных с двумя переменными:
\begin{equation*}
\frac{\partial^2 u(x,y)}{\partial x^2}+\frac{2\alpha}{x}\frac{\partial u(x,y)}{\partial x}=
\frac{\partial^2 u(x,y)}{\partial y^2}+\frac{2\beta}{y}\frac{\partial u(x,y)}{\partial y}
\end{equation*}
(обобщённое уравнение Эйлера--Пуассона--Дарбу или В--гиперболическое уравнение по терминологии И.А.~Киприянова) в открытой четверти плоскости $x>0,\  y>0$
при положительных параметрах $\beta>\alpha>0$ с краевыми условиями на осях координат (характеристиках)
\begin{equation*}
u(x,0)=f(x), \quad u(0,y)=g(y), \quad f(0)=g(0).
\end{equation*}
Предполагается, что решение $u(x,y)$ является непрерывно дифференцируемым в замкнутом первом квадранте, имеет непрерывные вторые производные в открытом квадранте, граничные функции $f(x), g(y)$ являются непрерывно дифференцируемыми.

Тогда, если решение поставленной задачи существует, то для него выполняются соотношения:
\begin{equation*}
\frac{\partial u}{\partial y}=0, y=0,\qquad  \frac{\partial u}{\partial x}=0, x=0,
\end{equation*}
\begin{equation*}
2^\beta \Gamma\left(\beta+\frac{1}{2}\right)\int\limits_0^1 f(xt)t^{\alpha+\beta+1}
(1-t^2)^{\frac{\beta -1}{2}}P_{-\alpha}^{1-\beta}(t)\,dt=
\end{equation*}
\begin{equation*}
=2^\alpha \Gamma\left(\alpha+\frac{1}{2}\right)\int\limits_0^1 g(xt)t^{\alpha+\beta+1}
(1-t^2)^{\frac{\alpha -1}{2}}P_{-\beta}^{1-\alpha}(t)\,dt,
\end{equation*}
\begin{equation*}
\Downarrow
\end{equation*}
\begin{equation*}
g(y)=\frac{2\Gamma(\beta+\frac{1}{2})}{\Gamma(\alpha+\frac{1}{2})
	\Gamma(\beta-\alpha)}y^{1-2\beta}
\int\limits_0^y x^{2\alpha-1}f(x)
(y^2-x^2)^{\beta-\alpha-1}x \,dx,
\end{equation*}
где  $P_\nu^\mu(z)$---функция Лежандра первого рода.

Перейдём к изложению результатов автора для ОП Бушмана--Эрдейи и их приложений к дифференциальным уравнениям с особенностями в коэффициентах.

Все рассмотрения ведутся ниже на полуоси. Поэтому будем обозначать через $L_2$ пространство $L_2(0, \infty)$ и $L_{2, k}$ весовое пространство $L_{2, k}(0, \infty)$.

Вначале распространим определения  (\ref{2BE1}--\ref{2BE2}) \  на  важный не исследованный ранее случай $\mu =1$.

\begin{dfn}
	Введём при $\mu =1$ операторы  Бушмана--Эрдейи нулевого порядка гладкости по формулам
	\begin{gather}
	\label{2BE01}
	(B_{0+}^{\nu,1}f)(x)={_1 S^{\nu}_{0+}f}=\frac{d}{dx}\int\limits_0^x P_\nu \left(\frac{x}{t}\right)f(t)\,dt,\\
	\label{2BE02}
	(E_{0+}^{\nu,1}f)(x)={_1 P^{\nu}_{-}}f=\int\limits_0^x P_\nu \left(\frac{t}{x}\right)\frac{df(t)}{dt}\,dt,\\
	\label{2BE03}
	(B_{-}^{\nu,1}f)(x)={_1 S^{\nu}_{-}}f=\int\limits_x^\infty P_\nu \left(\frac{t}{x}\right)\left(-\frac{df(t)}{dt}\right)\,dt,\\
	\label{2BE04}
	(E_{-}^{\nu,1}f)(x)={_1 P^{\nu}_{0+}}f=\left(-\frac{d}{dx}\right)\int\limits_x^\infty P_\nu \left(\frac{x}{t}\right)f(t)\,dt,
	\end{gather}
	где $P_\nu(z)=P_\nu^0(z)$---функция Лежандра.
\end{dfn}

Разумеется, при очевидных дополнительных условиях на функции в (\ref{2BE01})--(\ref{2BE04}) можно продифференцировать под знаком интеграла или проинтегрировать по частям.

\begin{theorem}\label{2fact1}
	Справедливы следующие формулы факторизации операторов Бушмана--Эрдейи на подходящих функциях через дробные интегралы Римана--Лиувилля и Бушмана--Эрдейи нулевого порядка гладкости:
	\begin{gather*}
	(B_{0+}^{\nu,\,\mu} f)(x)=I_{0+}^{1-\mu}~ {_1 S^{\nu}_{0+}f}, \quad (B_{-}^{\nu, \,\mu} f)(x)={_1 P^{\nu}_{-}}~ I_{-}^{1-\mu}f,\\
	(E_{0+}^{\nu,\,\mu} f)(x)={_1 P^{\nu}_{0+}}~I_{0+}^{1-\mu}f, \quad (E_{-}^{\nu, \, \mu} f)(x)= I_{-}^{1-\mu}~{_1 S^{\nu}_{-}}f.
	\end{gather*}
\end{theorem}

Эти важные формулы позволяют "разделить"\ параметры $\nu$ и $\mu$. Мы докажем, что операторы
\eqref{2BE01}--\eqref{2BE04} являются изоморфизмами пространств $L_2(0, \infty)$, если $\nu$
не равно некоторым исключительным значениям. Поэтому операторы \eqref{2BE1}--\eqref{2BE2}
по действию в пространствах типа $L_2$ в определённом смысле подобны операторам дробного интегродиффенцирования  $I^{1-\mu}$,
с которыми они совпадают при $\nu=0$. Далее операторы Бушмана--Эрдейи будут доопределены при всех значениях $\mu$.
Исходя из этого, введём следующее
\begin{dfn}
Число $\rho=1-{\rm Re}\,\mu $ назовём порядком гладкости операторов Бушмана--Эрдейи \eqref{2BE01}--\eqref{2BE04}.
\end{dfn}

Таким образом, при $\rho > 0$ (то есть при ${\rm Re}\,  \mu > 1$) операторы Бушмана--Эрдейи
являются сглаживающими, а при $\rho < 0$ (то есть при ${\rm Re}\, \, \mu < 1$) уменьшающими
гладкость в пространствах типа $L_2 (0, \infty)$. Операторы \eqref{2BE01}--\eqref{2BE04},
для которых $\rho = 0$, являются по данному определению операторами нулевого порядка гладкости.
Следует пояснить, что здесь под сглаживающими мы понимаем операторы, которые представимы в виде $A=D^k B$, где $k>0$,
а оператор $B$ ограничен в $L_2 (0, \infty)$.
Под уменьшающими гладкость при этом понимаются операторы, которые действуют из некоторого пространства $C^k(0, \infty), k>0$ в пространство Лебега $L_2 (0, \infty)$.

Важно отметить, что при некоторых специальных значениях параметров $\nu,~\mu$ операторы Бушмана--Эрдейи
сводятся к более простым. Так при значениях $\mu=-\nu$ или $\mu=\nu+2$ они являются
операторами Эрдейи--Кобера; при $\nu = 0$ операторами дробного интегродифференцирования
$I_{0+}^{1-\mu}$ или $I_{-}^{1-\mu}$; при $\nu=-\frac{1}{2}$, $\mu=0$ или $\mu=1$
ядра выражаются через эллиптические интегралы; при  $\mu=0$,  $x=1$, $v=it-\frac{1}{2}$  оператор $B_{-}^{\nu, \, 0}$
лишь на постоянную отличаются от преобразования Мелера--Фока. Таким образом, операторы Бушмана--Эрдейи первого рода являются обобщениями всех этих указанных классов стандартных интегральных операторов.

Будем рассматривать наряду с оператором Бесселя также тесно связанный с ним дифференциальный оператор
\begin{equation}
\label{275}
L_{\nu}=D^2-\frac{\nu(\nu+1)}{x^2}=
\end{equation}
\begin{equation*}
=\left(\frac{d}{dx}-\frac{\nu}{x}\right)
\left(\frac{d}{dx}+\frac{\nu}{x}\right)=\left(\frac{d}{dx}+\frac{\nu+1}{x}\right)
\left(\frac{d}{dx}-\frac{\nu+1}{x}\right),
\end{equation*}
который при $\nu \in \mathbb{N}$ является оператором углового момента из квантовой механики.
Их взаимосвязь устанавливают легко проверяемые формулы связи, приведём их.

Пусть пара ОП $X_\nu, Y_\nu$ сплетают $L_{\nu}$ и вторую  производную:
\begin{equation*}
X_\nu L_{\nu}=D^2 X_\nu,\qquad Y_\nu D^2 = L_{\nu} Y_\nu.
\end{equation*}
Введём новую пару ОП по формулам
\begin{equation*}
S_\nu=X_{\nu-1/2} x^{\nu+1/2},\qquad P_\nu=x^{-(\nu+1/2)} Y_{\nu-1/2},
\end{equation*}
тогда пара новых ОП $S_\nu, P_\nu$ сплетают оператор Бесселя и вторую производную:
\begin{equation*}
S_\nu B_\nu = D^2 S_\nu,\qquad P_\nu D^2 = B_\nu P_\nu.
\end{equation*}

Операторы нулевого порядка гладкости выделяются тем, что только для них можно доказать оценки в \textit{одном} пространстве типа $L_p(0,\infty)$.
При этом, учитывая структуру этих операторов, удобно пользоваться техникой преобразования Меллина и теоремой Слейтер (см. \cite{KS},\cite{SS}).

\begin{theorem} \label{2tmult}
	1. Операторы Бушмана--Эрдейи нулевого порядка гладкости действуют как умножение на мультипликатор в образах преобразования Меллина. Для их мультипликаторов справедливы формулы:
		\begin{gather*}
	m_{{_1S_{0+}^{\nu}}}(s)=\frac{\Gamma(-\frac{s}{2}+\frac{\nu}{2}+1) \Gamma(-\frac{s}{2}-\frac{\nu}{2}+\frac{1}{2})}{\Gamma(\frac{1}{2}-\frac{s}{2})\Gamma(1-\frac{s}{2})}=  \\
	=\frac{2^{-s}}{\sqrt{ \pi}} \frac{\Gamma(-\frac{s}{2}-\frac{\nu}{2}+\frac{1}{2}) \Gamma(-\frac{s}{2}+\frac{\nu}{2}+1)}{\Gamma(1-s)} ,\qquad {\rm Re}\, \, s < \min (2 + {\rm Re}\,  \,
	\nu, 1- {\rm Re}\, \, \nu); \\
	m_{{_1P_{0+}^{\nu}}}(s)=\frac{\Gamma(\frac{1}{2}-\frac{s}{2})\Gamma(1-\frac{s}{2})}{\Gamma(-\frac{s}{2}+\frac{\nu}{2}+1) \Gamma(-\frac{s}{2}-\frac{\nu}{2}+\frac{1}{2})},\qquad {\rm Re}\, \, s < 1; \\
	m_{{_1P_{-}^{\nu}}}(s)=\frac{\Gamma(\frac{s}{2}+\frac{\nu}{2}+1) \Gamma(\frac{s}{2}-\frac{\nu}{2})}{\Gamma(\frac{s}{2})\Gamma(\frac{s}{2}+\frac{1}{2})},\qquad {\rm Re}\,  \, s > \max({\rm Re}\,  \, \nu, -1-{\rm Re}\, \, \nu); \\
	m_{{_1S_{-}^{\nu}}}(s)=\frac{\Gamma(\frac{s}{2})\Gamma(\frac{s}{2}+\frac{1}{2})}{\Gamma(\frac{s}{2}+\frac{\nu}{2}+\frac{1}{2}) \Gamma(\frac{s}{2}-\frac{\nu}{2})},\qquad {\rm Re}\,  \, s >0.
	\end{gather*}

	2. Кроме того, выполняются следующие соотношения для мультипликаторов:
	\begin{gather*}
	m_{{_1P_{0+}^{\nu}}}(s)=1/m_{{_1S_{0+}^{\nu}}}(s),  \quad  m_{{_1P_{-}^{\nu}}}(s)=1/m_{{_1S_{-}^{\nu}}}(s),  \\
	m_{{_1P_{-}^{\nu}}}(s)=m_{{_1S_{0+}^{\nu}}}(1-s), \quad  m_{{_1P_{0+}^{\nu}}}(s)=m_{{_1S_{-}^{\nu}}}(1-s).
	\end{gather*}
	
	3. Справедливы следующие формулы для норм операторов Бушмана--Эрдейи нулевого порядка гладкости в $L_2$:
		\begin{gather*}
	\| _1{S_{0+}^{\nu}} \| = \| _1{P_{-}^{\nu}}\|= 1/ \min(1, \sqrt{1- \sin \pi \nu}),  \\
	\| _1{P_{0+}^{\nu}}\| = \| _1{S_{-}^{\nu}}\|= \max(1, \sqrt{1- \sin \pi \nu}).
	\end{gather*}

	4.  Нормы операторов \eqref{2BE01}--\eqref{2BE04} периодичны по $\nu$ с периодом 2, то есть $\|x^{\nu}\|=\|x^{\nu+2}\|$, где $x^{\nu}$ --- любой из операторов \eqref{2BE01}--\eqref{2BE04}.
	
	5. Нормы операторов ${_1 S_{0+}^{\nu}}$, ${_1 P_{-}^{\nu}}$ не ограничены в совокупности по $\nu$, каждая из этих норм не меньше $1$. Если $\sin \pi \nu \leq 0$, то эти нормы равны $1$. Указанные операторы неограничены в $L_2$ тогда и только тогда, когда $\sin \pi \nu = 1$ (или $\nu=(2k) + 1/2,~k \in \mathbb{Z}$).
	
	6. Нормы операторов ${_1 P_{0+}^{\nu}}$, ${_1 S_{-}^{\nu}}$
	ограничены в совокупности по $\nu$, каждая из этих норм не больше $\sqrt{2}$. Все эти операторы ограничены в $L_2$ при всех $\nu$. Если $\sin \pi \nu \geq 0$, то их $L_2$ -- норма равна 1. Максимальное значение нормы, равное $\sqrt 2 $, достигается тогда и только тогда, когда $\sin \pi \nu = -1$ (или $\nu= -1/2+(2k) ,~k \in \mathbb{Z}$).
\end{theorem}

Важнейшим свойством операторов Бушмана--Эрдейи нулевого порядка гладкости является их унитарность при целых $\nu$.
Отметим, что при интерпретации $L_{\nu}$ как оператора углового момента в квантовой механике, параметр $\nu$ как раз и принимает целые неотрицательные значения.
Сформулируем один из основных результатов данной главы.

\begin{theorem}\label{2tunit}
Для унитарности в $L_2$ операторов \eqref{2BE01} -- \eqref{2BE04} необходимо и достаточно, чтобы число $\nu$ было целым.
В этом случае пары операторов
	$({_1 S_{0+}^{\nu}}$, ${_1 P_{-}^{\nu}})$ и  $({_1 S_{-}^{\nu}}$, ${_1 P_{0+}^{\nu}})$
	взаимно обратны.
\end{theorem}

Перейдём к построению операторов преобразования, унитарных при  всех $\nu$. Такие операторы определяются по формулам:
\begin{gather}
S_U^{\nu} f = - \sin \frac{\pi \nu}{2}\  {_2S^{\nu}}f+ \cos \frac{\pi \nu}{2}\  {_1S_-^{\nu}}f, \label{2.614} \\
P_U^{\nu} f = - \sin \frac{\pi \nu}{2}\  {_2P^{\nu}}f+ \cos \frac{\pi \nu}{2}\  {_1P_-^{\nu}}f. \label{2.615}
\end{gather}
Для любых значений $\nu \in \mathbb{R}$ они являются линейными комбинациями операторов преобразования Бушмана--Эрдейи 1 и 2 рода нулевого порядка гладкости.
Их можно назвать операторами Бушмана--Эрдейи третьего рода. В интегральной форме эти операторы имеют вид:

\begin{gather}
S_U^{\nu} f = \cos \frac{\pi \nu}{2} \left(- \frac{d}{dx} \right) \int\limits_x^{\infty} P_{\nu}\left(\frac{x}{y}\right) f(y)\,dy +  \label{2.616}\\
+ \frac{2}{\pi} \sin \frac{\pi \nu}{2} \left(\int\limits_0^x (x^2-y^2)^{-\frac{1}{2}}Q_{\nu}^1 \left(\frac{x}{y}\right) f(y)\,dy  \right.-\nonumber\\
 -  \int\limits_x^{\infty} (y^2-x^2)^{-\frac{1}{2}}\mathbb{Q}_{\nu}^1 \left(\frac{x}{y}\right) f(y)\,dy \Biggl. \Biggr), \nonumber \\
 P_U^{\nu} f = \cos \frac{\pi \nu}{2}  \int\limits_0^{x} P_{\nu}\left(\frac{y}{x}\right) \left( \frac{d}{dy} \right) f(y)\,dy - \label{2.617} \\
  -\frac{2}{\pi} \sin \frac{\pi \nu}{2} \left( - \int\limits_0^x (x^2-y^2)^{-\frac{1}{2}}\mathbb{Q}_{\nu}^1\left(\frac{y}{x}\right) f(y)\,dy   \right.-\nonumber\\
- \int\limits_x^{\infty} (y^2-x^2)^{-\frac{1}{2}} Q_{\nu}^1 \left(\frac{y}{x}\right) f(y)\,dy \Biggl. \Biggr). \nonumber
\end{gather}

\begin{theorem}\label{2unit} Операторы \eqref{2.614}--\eqref{2.615} или  \eqref{2.616}--\eqref{2.617} при всех $\nu$ являются унитарными, взаимно сопряжёнными и обратными в $L_2$. Они являются
	сплетающими и действуют по формулам \eqref{275}. При этом $S_U^{\nu}$ является оператором типа Сонина (Сонина--Катрахова), а $P_U^{\nu}$ --- типа Пуассона (Пуассона--Катрахова).
\end{theorem}

\section{Дробные степени операторов Бесселя}

В этом пункте мы рассмотрим второй важный класс обобщённых операторов дробного интегродифференцирования --- это дробные степени операторов Бесселя.

Рассмотрим вещественные степени сингулярного дифференциального оператора Бесселя
\begin{equation}\label{Bess}
B_\nu= D^2+\frac{\nu}{x}D,\qquad \nu\geq 0
\end{equation}
на вещественной полуоси $(0,\infty)$.

\begin{dfn}
Пусть $\alpha>0$, $f(x)\in C^{[2\alpha]+1}(0,\infty)$. Дробную степень оператора Бесселя на полуоси $(0,\infty)$ определим, следуя работам  \cite{Kosovo}, \cite{Sita1}, \cite{Sita2}, \cite{Ida} формулой
\begin{equation*}
(IB_{\nu,-}^{\alpha}\,f)(x){=}
\end{equation*}
\begin{equation}
\label{Bess3}
=\frac{1}{\Gamma(2\alpha)}\int\limits_x^{+\infty}\left(\frac{y^2-x^2}{2y}
\right)^{2\alpha-1}\,_2F_1\left(\alpha+\frac{\nu-1}{2},\alpha;2\alpha;1-\frac{x^2}{y^2}\right)f(y)dy.
\end{equation}
Для краткости будем также называть выражение \eqref{Bess3} \textbf{дробным интегралом Бесселя на полуоси}.
\end{dfn}

В работе \cite{McBride}  введены пространства, приспособленные для работы с операторами вида \eqref{Bess3}:
\begin{equation*}
F_p=\left\{\varphi\in C^\infty(0,\infty):x^k\frac{d^k\varphi}{dx^k}\in L^p(0,\infty)\,\, {\text{для}}\,k=0,1,2,...\right\},\qquad 1\leq p<\infty,
\end{equation*}
\begin{equation*}
F_\infty=\left\{\varphi\in C^\infty(0,\infty):x^k\frac{d^k\varphi}{dx^k}\rightarrow0 \,\, {\text{при}}\, x\rightarrow0+
 {\text{,\,\,и при}}\, x\rightarrow\infty\,{\text{для}}\,k=0,1,2,...\right\}
\end{equation*}
и
\begin{equation*}
F_{p,\mu}=\left\{\varphi: x^{-\mu}\varphi(x)\in F_p\right\},\qquad 1\leq p\leq \infty,\qquad \mu\in\mathbb{C}.
\end{equation*}
Кроме того, в \cite{McBride}  доказано, что \eqref{Bess3} имеет обратный оператор.

\begin{dfn}
 \textit{Дробную производную Бесселя на полуоси} определим равенством
\begin{equation}
\label{DrobessDer}
(DB_{\nu,-}^\alpha f)(x)=B_\nu^n(IB_{\nu,-}^{n-\alpha}f)(x),\qquad \alpha>0.
\end{equation}
\end{dfn}

Можно показать, что на подходящем классе функций $DB_{\nu,-}^\alpha$ есть левый обратный оператор к $IB_{\nu,-}^\alpha$.

 Используя формулу, связывающую гипергеометрическую функцию Гаусса и функцию Лежандра вида
\begin{equation*}
_2F_1(a,b;2b;z)=2^{2b-1}\Gamma\left(b+\frac{1}{2}\right)\, z^{\frac{1}{2}-b}(1-z)^{\frac{1}{2}\left(b-a-\frac{1}{2}\right)}P_{a-b-\frac{1}{2}}^{\frac{1}{2}-b}
\left[\left(1-\frac{z}{2}\right)\sqrt{1-z}\right]
\end{equation*}
(см. формулу 15.4.8 на стр. 561 из \cite{Abramowitz}), мы получим
\begin{equation*}
\,_2F_1\left(\alpha+\frac{\nu-1}{2},\alpha;2\alpha;1-\frac{x^2}{y^2}\right)=
\end{equation*}
\begin{equation*}
=2^{2\alpha-1}\Gamma\left(\alpha+\frac{1}{2}\right)\, \left(\frac{y^2-x^2}{y^2}\right)^{\frac{1}{2}-\alpha}\left(\frac{y}{x}
\right)^{\frac{\nu}{2}}P_{\frac{\nu}{2}-1}^{\frac{1}{2}-\alpha}
\left[\frac{1}{2}\left(\frac{x}{y}+\frac{y}{x}\right)\right],
\end{equation*}
и сможем записать \eqref{Bess3} в виде
\begin{equation*}
(B_{\nu,-}^{-\alpha}f)(x)
=\frac{\Gamma\left(\alpha+\frac{1}{2}\right)}{\Gamma(2\alpha)}\int\limits_x^b(y^2-x^2)^{\alpha-\frac{1}{2}}
\,\left(\frac{y}{x}\right)^{\frac{\nu}{2}}P_{\frac{\nu}{2}-1}^{\frac{1}{2}-\alpha}\left[\frac{1}{2}
\left(\frac{x}{y}+\frac{y}{x}\right)\right]f(y)dy.
\end{equation*}

Выражение дробных интегралов Бесселя через функции Лежандра является полезным и является упрощением первоначального определения, так как гипергеометрическая функция Гаусса зависит от трёх параметров,
а функция Лежандра --- от двух.

Далее изучается дробный интеграл Бесселя в форме \eqref{Bess3}, который в частном случае соответствует дробному интегралу Лиувилля.
Существует также версия дробного интеграла Бесселя на конечном отрезке с интегрированием по промежутку $(0,x)$, которая в частном случае является дробным интегралом Римана--Лиувилля,
а также их дальнейшие модификации, см. \cite{Kosovo}, \cite{Sita1}, \cite{Sita2}.

Определение \eqref{Bess3} дано выше при ограничениях на функцию, близким к оптимальным в указанном классе, однако далее мы будем для простоты предполагать, что рассматриваются бесконечно дифференцируемые функции, финитные на полуоси, то есть их носитель  есть ${\rm supp}\,{f(x)}=[a,b], 0<a<b<\infty$.

Приведём основные свойства оператора \eqref{Bess3}, во-первых, демонстрирующие связь дробного интеграла Бесселя на полуоси с дробным интегралом Лиувилля и с дробным интегралом Сайго, во-вторых, показывающие, что при дополнительных условиях оператор \eqref{Bess3} при $\alpha=1$ обращает оператор Бесселя \eqref{Bess}, и, наконец, найдем дробный интеграл Бесселя на полуоси \eqref{Bess3} от степенной функции.

\begin{utv} При $\nu=0$ дробный интеграл Бесселя на полуоси $B_{0,-}^{-\alpha}$ сводится к  дробному интегралу Лиувилля, определённому формулой (5.3) стр. 85 из \cite{SKM}, а именно, справедлива формула
\begin{equation*}
(IB_{0,-}^{\alpha}f)(x)=\frac{1}{\Gamma(2\alpha)}\int\limits_x^{\infty}(y-x)^{2\alpha-1}f(y)dy=(I_{-}^{2\alpha}f)(x).
\end{equation*}
\end{utv}

\begin{utv} Имеет место равенство
\begin{equation*}
(IB_{\nu,-}^{\alpha}f)(x)=\frac{1}{2^{2\alpha}}J_{x^2}^{2\alpha,\frac{\nu{-}1}{2}-\alpha,-\alpha}\left(x^{\frac{\nu{-}1}{2}}f(\sqrt{x})\right),
\end{equation*}
где
\begin{equation}
\label{Saigo1}
J_{x}\,^{\gamma,\beta,\eta}f(x)=\frac{1}{\Gamma(\gamma)}\int\limits_x^\infty(t-x)^{\gamma-1}t^{-\gamma-\beta}\,_2F_1\left(\gamma+\beta,-\eta;\gamma;1-\frac{x}{t}\right)f(t)dt,
\end{equation}
--- дробный интеграл Сайго   (см. \cite{Saigo}, \cite{Repin}). В \eqref{Saigo1} $\gamma>0,\beta,\theta$ --- вещественные числа.
\end{utv}
\begin{utv} При $\lim\limits_{x\rightarrow +\infty}g(x)=0$, $\lim\limits_{x\rightarrow +\infty}g'(x)=0$ получим, что
\begin{equation*}
(IB_{\nu,-}^{-1}B_\nu g)(x)=g(x).
\end{equation*}
\end{utv}

\begin{utv} При $x>0$ и $m+2\alpha+\nu<1$ справедлива формула
\begin{equation}
\label{Prop4}
IB_{\nu,-}^{\alpha}\,x^m=2^{-2\alpha}\,\Gamma\left[
                                                           \begin{array}{cc}
                                                             -\alpha-\frac{m}{2}, & -\frac{\nu-1}{2}-\alpha-\frac{m}{2} \\
                                                             \frac{1-\nu-m}{2}, & -\frac{m}{2}  \\
                                                           \end{array}
                                                         \right]\, x^{2\alpha+m}.
\end{equation}
\end{utv}

Отметим важность полученной формулы, устанавливающей, что дробный интеграл Бесселя переводит одну степенную функцию в другую, так как это позволяет распространить его определение на произвольные степенные ряды.
Явный вид константы в формуле \eqref{Prop4} в форме отношения гамма--функций показывает, что дробный интеграл Бесселя является оператором дробного дифференцирования типа Гельфонда--Леонтьева \cite{SKM}.

Теперь выведем  формулу преобразования Меллина от дробного интеграла Бесселя на полуоси $B_{\nu,-}^{-\alpha}$.

Преобразование Меллина функции $f$ определяется формулой
 \begin{equation*}
 Mf(s)=f^*(s)=\int\limits_0^\infty x^{s-1}f(x)dx.
 \end{equation*}

\begin{theorem} Пусть $\alpha>0$. Преобразования Меллина от дробного интеграла и дробной производной Бесселя на полуоси имеют вид
\begin{equation*}
  ((IB_{\nu,-}^{\alpha}f)(x))^*(s)=\frac{1}{2^{2\alpha}}\,\,\Gamma\left[\begin{array}{cc}
\frac{s}{2}, & \frac{s}{2}-\frac{\nu-1}{2} \\
\alpha+\frac{s}{2}-\frac{\nu-1}{2}, & \alpha+\frac{s}{2}  \\
                                                           \end{array}
                                                         \right] f^*(2\alpha+s),
\end{equation*}
\begin{equation*}
  ((DB_{\nu,-}^{\alpha}f)(x))^*(s)=2^{2\alpha}\,\Gamma\left[\begin{array}{cc}
\frac{s}{2}, & \frac{s}{2}-\frac{\nu-1}{2} \\
\frac{s}{2}-\alpha-\frac{\nu-1}{2}, & \frac{s}{2}-\alpha  \\
\end{array} \right] f^*(s-2\alpha).
\end{equation*}
\end{theorem}

\begin{theorem} Для дробного интеграла и дробной производной Бесселя на полуоси при $\alpha,\beta>0$ справедливо полугрупповое свойство
\begin{equation*}
    IB_{\nu,-}^{-\alpha}IB_{\nu,-}^{\beta}f= IB_{\nu,-}^{\alpha+\beta}f,
\end{equation*}
\begin{equation*}
    DB_{\nu,-}^{\alpha}DB_{\nu,-}^{\beta}f= DB_{\nu,-}^{\alpha+\beta}f.
\end{equation*}
\end{theorem}

Отметим, что в литературе по интегральным преобразованиям с параметром часто полугрупповое свойство называется индексным законом.

В заключение отметим, что конструкции дробных степеней оператора Бесселя вида \eqref{Bess3} и \eqref{DrobessDer} могут быть применены к исследованию различных дифференциальных уравнений дробного порядка,
а также в теории операторов преобразования, см. \cite{Sita3}--\cite{Sita10}, \cite{SS1}, \cite{SS2}.
Дифференциальное уравнение с  дробной степенью оператора Бесселя используется при моделировании случайного блуждания частицы (см. \cite{Garra1}, \cite{Garra2}).
Другой подход к определению дробных степеней оператора Бесселя см. в \cite{McBride}.
Исследованию различных степеней оператора гиперболического типа, содержащего операторы Бесселя, посвящены работы \cite{Sh1}--\cite{Sh4}.

\end{document}